# Spectral, entropy and bifurcation analysis of the dynamics of a reverse-flow system

*Marek Berezowski, Natalia Koziol, Marcin Lawnik*

Cracow University of Technology, Faculty of Chemical Engineering and Technology, 30-155 Cracow, ul. Warszawska 24, Poland

Corresponding author: marek.berezowski@pk.edu.pl

**Abstract**

The work concerns the spectral, entropy and bifurcation analysis of the dynamics of a reverse-flow system. The existence of chaotic oscillations was demonstrated in a wide range of changes in the parameters of the model.

The model of such a system can be compared to a periodically stimulated pendulum. The periodic apparatus here is a tubular reactor, which can generate its own oscillations. Reverse flow is cyclical stimulation. As a result, both of these devices can start to oscillate chaotically, and they actually do. The purpose of using reverse flow is to increase the conversion degree. This is because after reverse, the cold stream of raw material falls into the heated reactor inlet.

## 1. Introduction

The dynamics of chemical reactors was widely discussed in [1-17], where it was proved that in the steady state the concentration and temperature of the flux may oscillate in a manner that is more or less complex. The oscillations may have a multiperiodic, quasiperiodic or chaotic nature, depending on the values of reactor parameters and the concentration and temperature of the reacting flux at the initial state.

The dynamic behaviour of periodically forced chemical reactors (e.g. with cyclically permuting feed and discharge positions or changing the feed flow direction) was studied in [1,4-5,7-9,12-17, 23]. It was shown that in these networks the type of oscillations also depends on the switching time (i.e. the time interval occurring between two successive permutations of the reactor order or flow reversals, respectively).

This paper is focused on the analysis of the dynamics of a system consisting of catalyst (psudohomogeneous), non-adiabatic tubular chemical reactor with longitudinal dispersion. An interesting work on the axial dispersion tubular reactor can be found here: [20]. The reactor is fed with flux of raw materials with cyclically reversing the flow directions. The calculations results indicated that different types of temperature and concentration oscillations may occur in the analyzed system, including chaotic oscillations. The analysis was based on bifurcation diagrams, amplitude spectrum and informatic entropy.

The research carried out in this paper is of a general nature. The presented methods and conclusions can and should be applied at the stage of designing the reactor system. It is particularly important that chaotic oscillations, which are usually unfavorable in terms of processes, may occur in the tested system. Therefore, the operating parameters of the system should be selected in such a way as to eliminate these oscillations in practice. We learn about the degree of danger that can be caused by all kinds of oscillations, not only chaotic ones, from the spectral analysis based on the frequencies and amplitudes of individual harmonics. In turn, the entropy analysis allows the identification of the physical system under study. This is necessary at the system modeling stage. Bifurcation analysis, on the other hand, gives information about the nature and values of the system state variables, in this case about the conversion degree and temperature.

## 2. The reactor model

The analysis concerns of a pseudohomogeneous, tubular, non-adiabatic chemical reactor with longitudinal dispersion. The balance equations of the apparatus are as follows:

mass balance:

$$\frac{\partial \alpha(\xi,\tau)}{\partial \tau} + \frac{\partial \alpha(\xi,\tau)}{\partial \xi} = \frac{1}{Pe_M}\frac{\partial^2 \alpha(\xi,\tau)}{\partial \xi^2} + \Phi_1 \qquad (1)$$

heat balance:

$$Le\frac{\partial \Theta(\xi,\tau)}{\partial \tau} + \frac{\partial \Theta(\xi,\tau)}{\partial \xi} = \frac{1}{Pe_H}\frac{\partial^2 \Theta(\xi,\tau)}{\partial \xi^2} + \Phi_2 \qquad (2)$$

where the dimensionless position in the reactor is in the range:

$$0 \leq \xi \leq 1. \qquad (3)$$

$\alpha(\xi,\tau)$ is the degree of conversion along the position $\xi$ at any time $\tau$, $\Theta(\xi, \tau)$ is the dimensionless temperature along the position $\xi$ at any time $\tau$. The variables $\alpha$ and $\Theta$ are related to the concentration and temperature in the raw material stream (see *Notation*).

Assuming that the *n*-order A→B type reaction takes place in the reactor, the kinetic and heat exhange functions are as follows:

$$\Phi_1 = Da(1-\alpha)^m \exp\left(\gamma\frac{\beta\Theta}{1+\beta\Theta}\right) \qquad (4)$$

$$\Phi_2 = \Phi_1(\alpha,\Theta) + \delta(\Theta_H - \Theta). \qquad (5)$$

The boundary conditions prescribed by the balance equations are of the Danckwerts type:

$$\alpha(0,\tau) = \frac{1}{Pe_M}\frac{d\alpha(0,\tau)}{d\xi} \qquad (6)$$

$$\frac{d\alpha(1,\tau)}{d\xi} = 0 \qquad (7)$$

$$\Theta(0,\tau) = \frac{1}{Pe_H}\frac{d\Theta(0,\tau)}{d\xi} \qquad (8)$$

$$\frac{d\Theta(1,\tau)}{d\xi} = 0. \tag{9}$$

In the above model, cyclic transfer of the raw material stream (reverse flow) was used. The product is collected at the reactor outlet according to the following relations:

$$\alpha_{out} = IO\alpha(0,\tau)+(1\text{-}IO)\alpha(1,\tau) \tag{10}$$

$$\Theta_{out} = IO\Theta(0,\tau)+(1\text{-}IO)\Theta(1,\tau) \tag{11}$$

where the *IO* index specifies the current flow direction of the reacting stream and has a value of 0 or 1. This direction changes cyclically at times $\tau=k\tau_r$, where $k = 1,2,3,$ ... and $\tau_r$ is the switching time.

The following parameter values were adopted for the numerical calculations carried out in this paper: γ=15, β=2, m=1.5, δ=3, $\Theta_H = 0$, $Pe_M = 50$, $Pe_H = 50$, *Le*=1.

## 3. Bifurcation analysis

The shift flow of the raw material stream at regular intervals $\tau_r$ causes that the tested system can be treated as a discrete system. This means that in order to evaluate the system's behavior it is sufficient to observe the $\alpha_{out}$ or $\Theta_{out}$ variables at the moments of the transfer. The analysis carried out as part of this study was performed with the use of tools for discrete systems. Depending on the values of the parameters, the above variables behave in a stationary or oscillatory manner. These oscillations can be periodic or non-periodic, including chaotic. Fig. 1 presents a bifurcation diagram of the state variable $\alpha_{out}$ depending on the switching time $\tau_r$.

It is formed as follows. For a given value of $\tau_r$, a discrete time series is generated. Then all values of the $\alpha_{out}$ variable read from the time series are plotted on the diagram. If in the time series this variable has the same value, then in the diagram we see one point (e.g. for $\tau_r= 4$, $\tau_r= 13.5$ in Fig. 1). This is called stationary point. If a variable has *N* different values in the time series, then we see *N* points in the diagram. These are *N*-period oscillations. If $N = \infty$, then we are dealing with chaos or quasiperiodic oscillations. The other bifurcation diagrams are created in the same way. Only the so-called the bifurcation parameter, which is *Da* in Fig. 2.

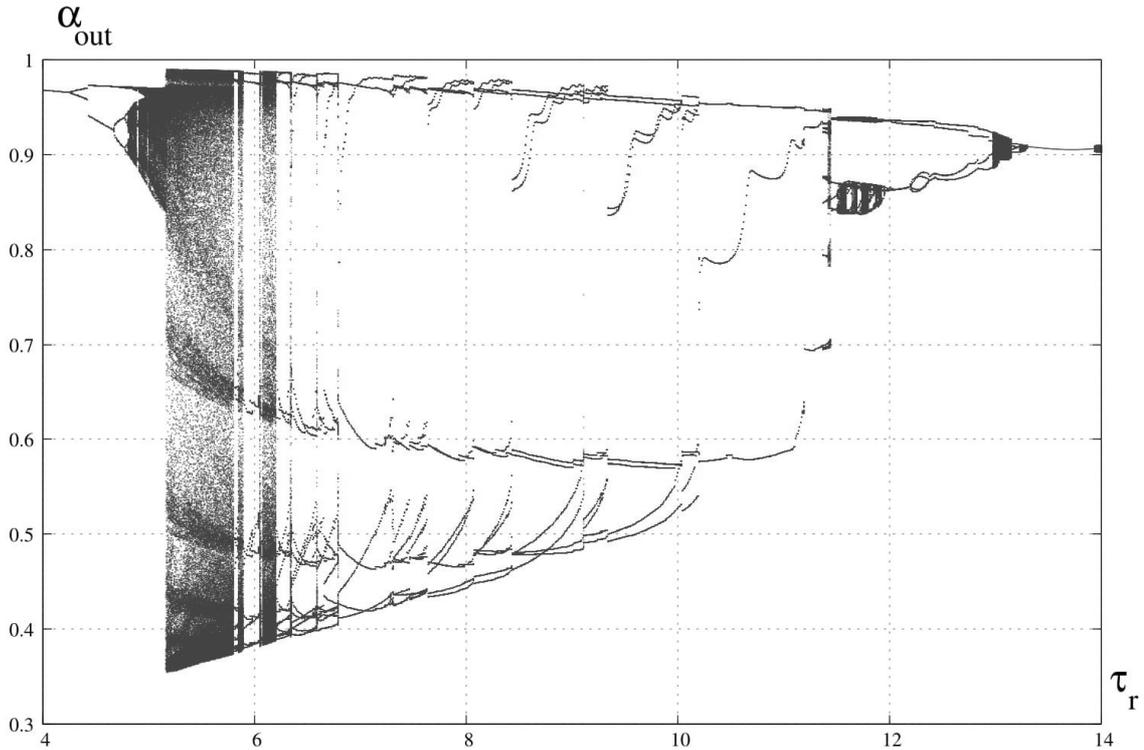

Fig. 1. Bifurcation diagram for *Da* = 0.13

In the steady state, the time series is characterized by the same value of the variable and we see a single point in the diagram. In the case of periodic oscillations, the diagram shows the number of points corresponding to the periodicity of a given discrete waveform (a given periodic orbit). In the case of chaotic oscillations, an infinite number of points should appear in the diagram. Fig. 2 presents a bifurcation diagram of the variable $\alpha_{out}$ depending on the number of *Da*.

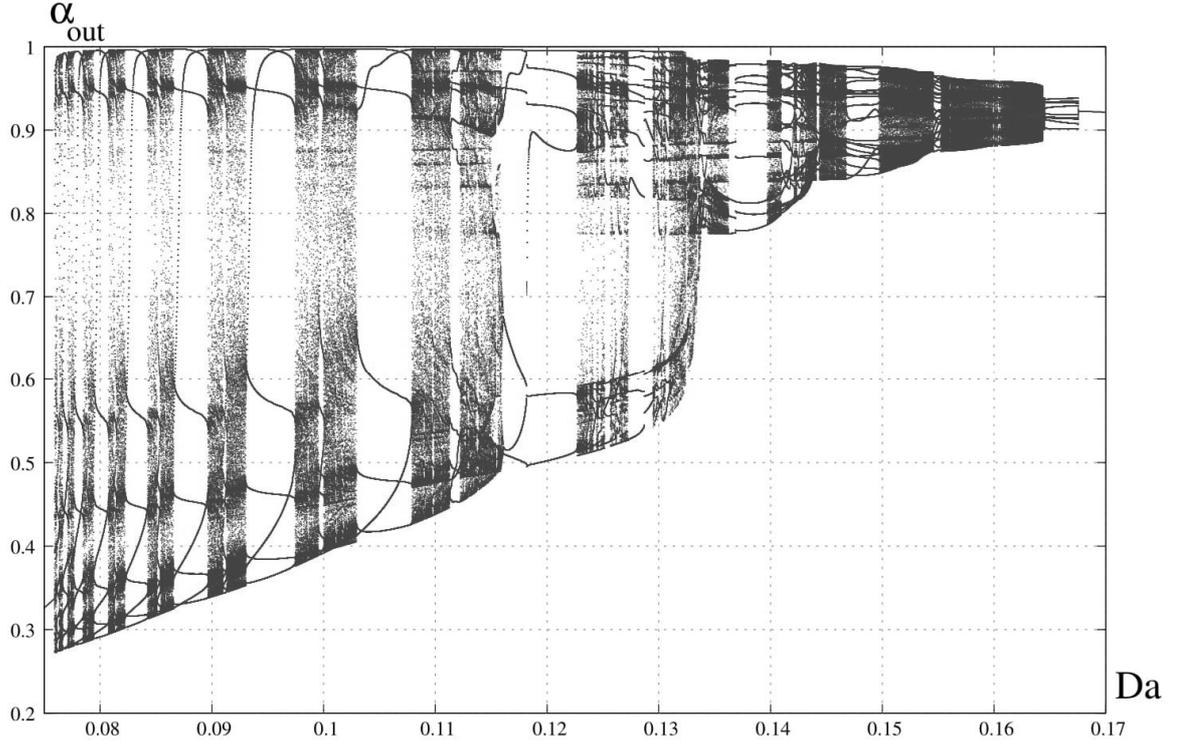

Fig. 2. Bifurcation diagram for $\tau_r = 5.5$

## 4. Spectral analysis

As part of this study, the amplitude spectrum of the degree of conversion $\alpha_{out}$ was analyzed according to the relationship:

$$X[k] = \frac{1}{N}\sum_{n=0}^{N-1} \alpha_{out}(n) exp\left(-\frac{2\pi i k n}{N}\right) \qquad (12)$$

where $\alpha_{out}(n)$ is the degree of conversion in the product stream in the nth step of the transfer, X [k] is the kth harmonic and N is the number of samples. Similarly to the bifurcation analysis, also here a diagram was determined, in which harmonics as a function of switching time $\tau_r$ were marked (Fig. 3) [18]. The principle of creating this diagram is the same as creating a bifurcation diagram, except that it does not use the time series, but the sequence of harmonics generated by the formula (12).

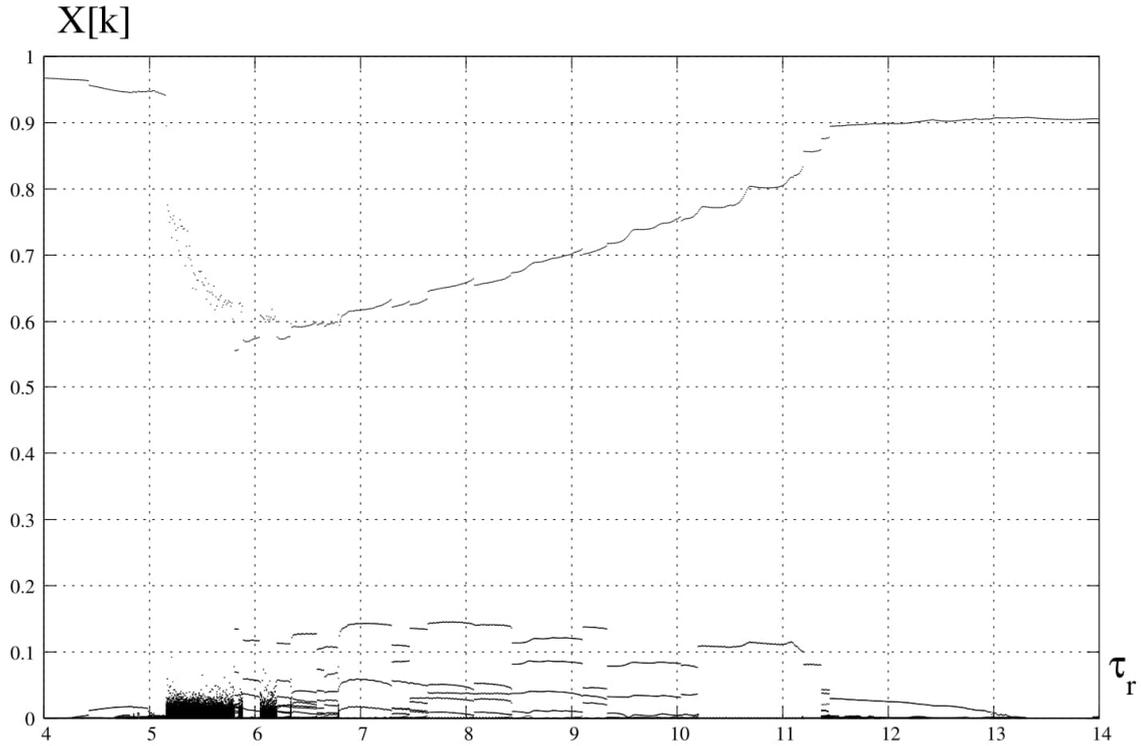

Fig. 3. Diagram of the harmonics of the amplitude spectrum for Da = 0.13

From this diagram it can be concluded that chaotic solutions occur for $4.8 < \tau_r < 6.2$. We draw a similar conclusion from the analysis of the diagram in Fig. 1, where the same chaotic interval is visible. It is obvious. If these ranges do not coincide, it would mean that the calculations are incorrect.

Fig. 4 shows the harmonic distribution for $\tau_r = 5.5$. This distribution applies to a chaotic orbit. Theoretically, there should be infinitely many harmonics in this case.

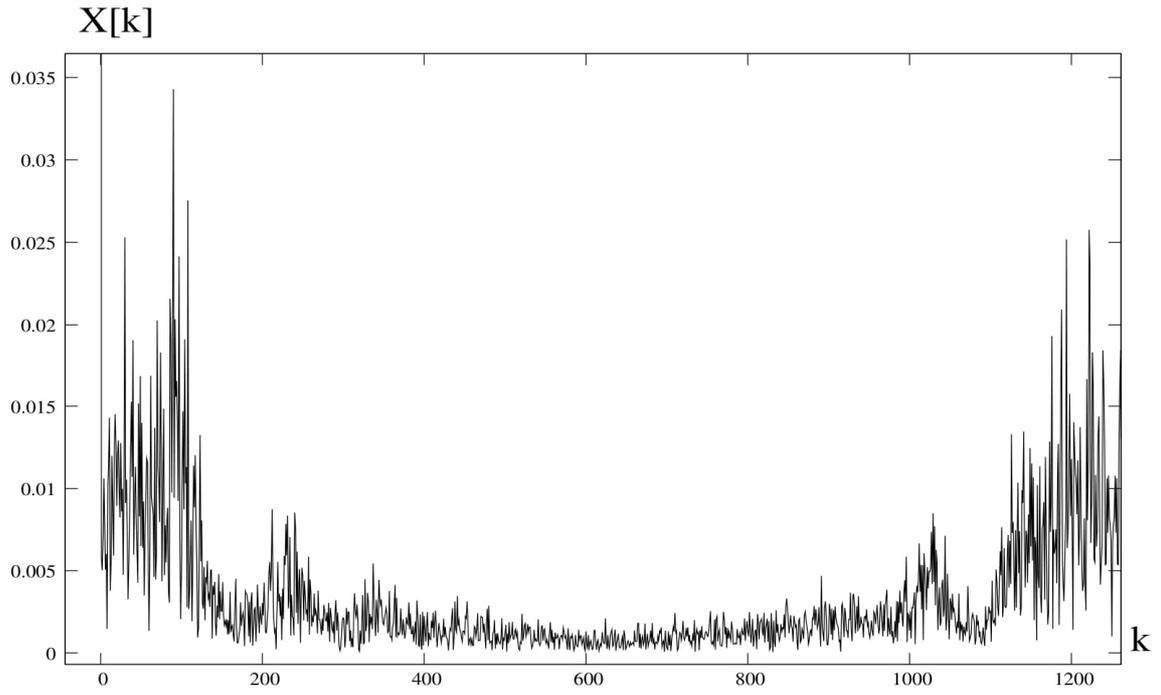

Fig. 4. Harmonics of the amplitude spectrum for Da = 0.13 and $\tau_r = 5.5$

Fig. 5 shows the harmonic distribution for $\tau_r = 6.5$. This distribution is for a stable periodic orbit.

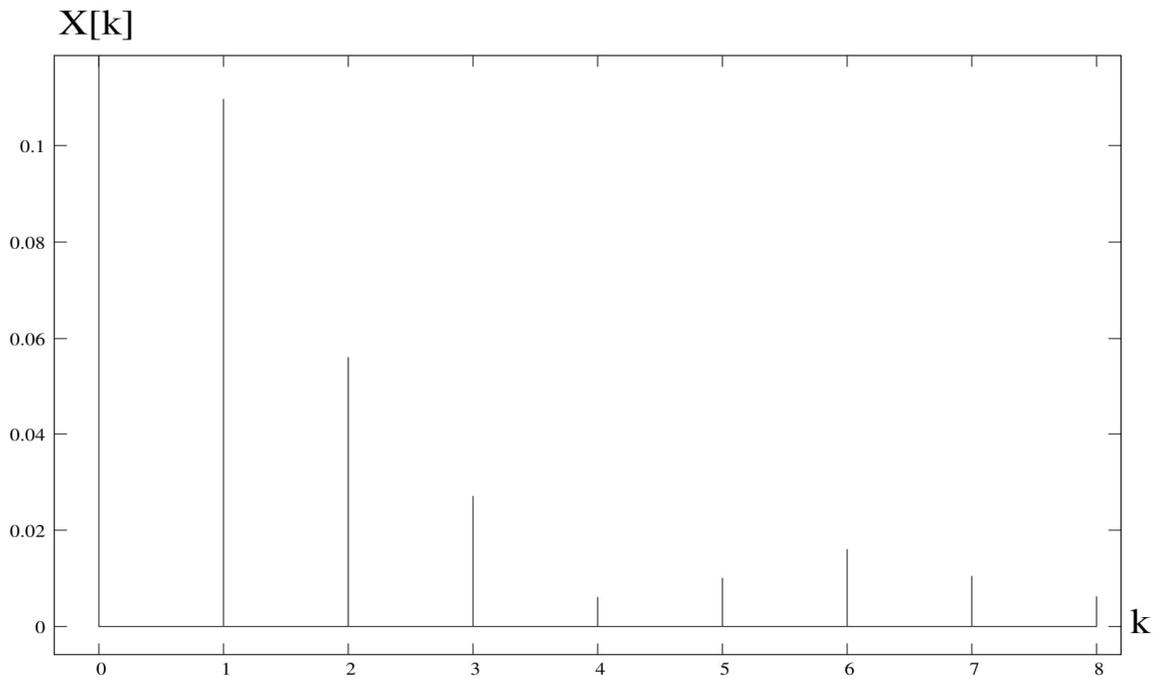

Fig. 5. Harmonics of the amplitude spectrum for $Da = 0.13$ and $\tau_r = 6.5$

(In Figs. 4 and 5, high values of zero harmonics have been limited due to the legibility of the graphs).

Additionally, Fig. 6 shows the cross-section of the Poincare chaotic attractor for $\tau_r = 5.5$. Due to the discrete nature of the system, this cross-section is of the Henon attractor type.

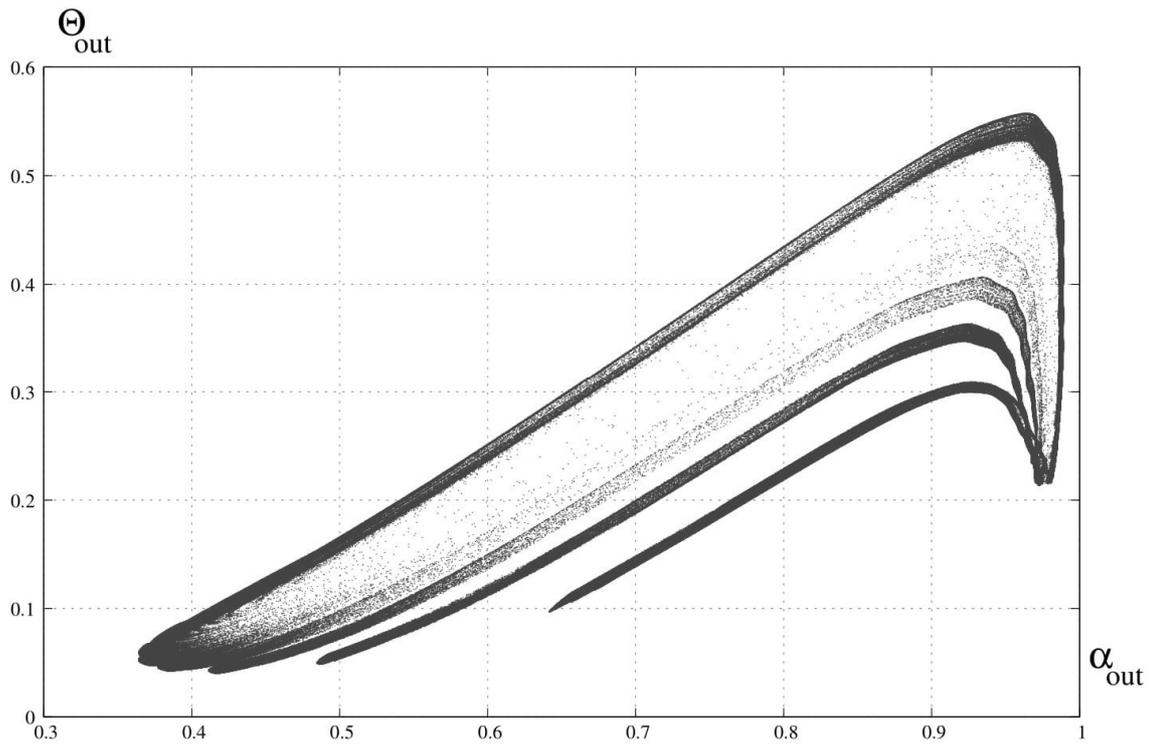

Fig. 6. Poincare cross-section for Da = 0.13 and $\tau_r = 5.5$

Numerical calculations showed that the tested system is characterized by the so-called a multiple of dynamic states. Each of these states depends on the initial conditions in the reactor. Fig. 7 shows another bifurcation diagram, which was created for different initial conditions than the diagram from Fig. 1.

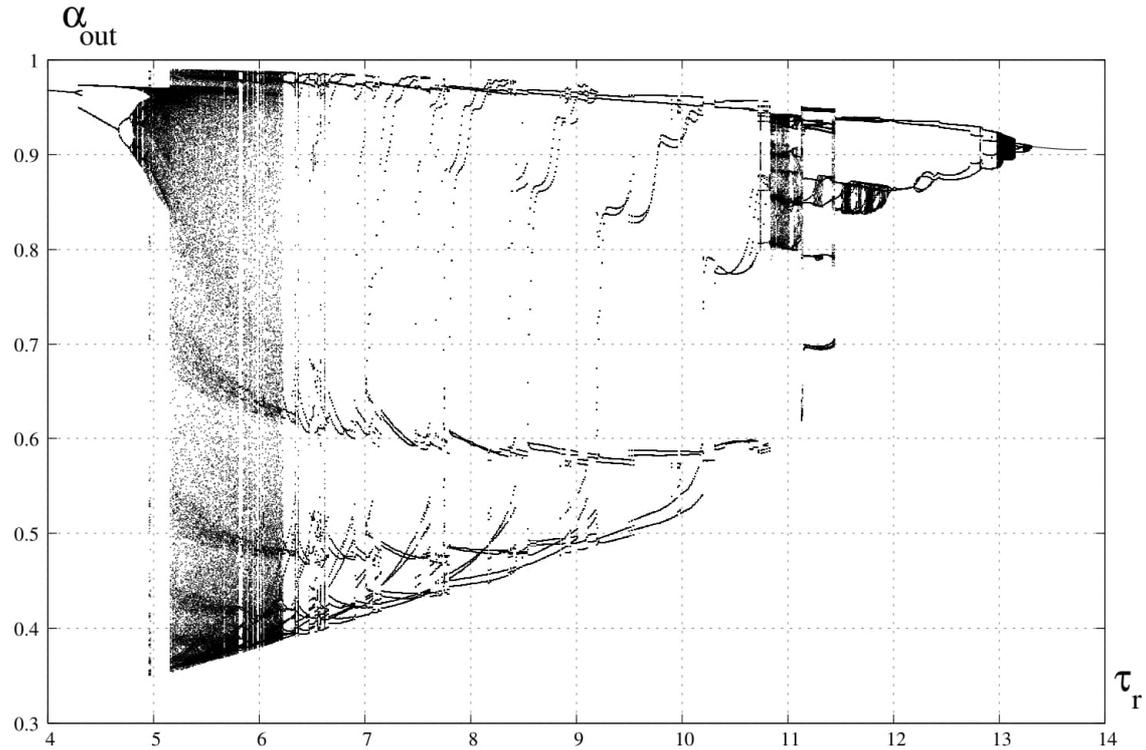

Fig. 7. Bifurcation diagram for *Da*=0.13

Comparing the two diagrams, it can be seen that in the latter case, chaos also occurs for $\tau_r = 6.5$ and for $\tau_r = 11$. In the first case, for $\tau_r = 6.5$ and for $\tau_r = 11$, the system generated periodic orbits. It should be remembered that the values of all system parameters in both cases are the same. In the first case (Fig. 1), the diagram was obtained assuming the following initial values: α(ξ,0)=0.9, Θ(ξ,0)=0.2. In the second case (Fig. 7), α(ξ,0)=0.2, Θ(ξ,0)=0.1 was assumed. How sensitive the tested system is to the initial values of the variables is shown in Fig. 8. It shows a map of the initial conditions determining the periodic solutions of the reactor (dots). Empty spaces determine chaotic solutions. By initial condition is meant here the state of the reactor at time zero along the entire length of the tube.

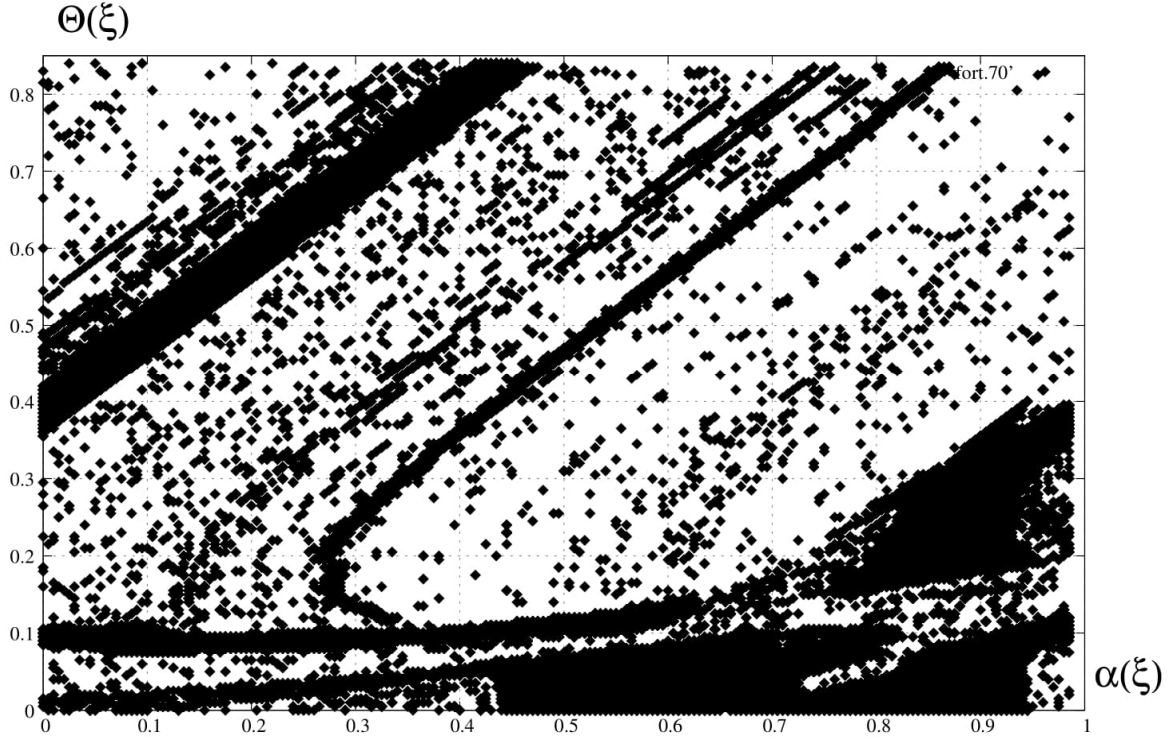

Fig. 8. Map of periodic and chaotic solutions

## 5. Entropy analysis

The last method of testing the discussed reactor is the entropy analysis [19]. It allows to define the range of parameters for which the system generates the most information about itself. This can be useful in the reactor modeling stage. The concept of entropy is used here to define the amount of information contained in the system [21]. The amount of entropy in a probability distribution is known as the information content of this distribution. Information is defined in terms of the likelihood of certain events occurring in the past; the greater the earlier uncertainty of such an event, the greater the amount of information obtained when such an event occurs, which suggests that the measure would range from zero to infinity [21]. Shannon [22] was the first to derive a general formula for measuring entropy in any set of probabilities. It is written as:

$$E = -\sum_{i=1}^{N} p_i \log_2 p_i. \qquad (13)$$

This formula is an indicator of the amount of information about a given system generated by this system. This is called information entropy, where N is the number of samples, and $p_i$ is the probability of the *i*-th number in the set of samples. For stationary case $E = 0$, for 1-period orbit $E = 1$, for 2-period orbit $E = 2$, for 3-period orbit $E = 1.585$, etc. General:

$$E = \log_2 M \qquad (14)$$

where *M* is the periodicity of a given orbit. It follows that in the case of a chaotic orbit $E=\infty$. In practice, these values may be slightly different due to the accuracy of numerical calculations. It is obvious that with low computational accuracy, more than one number may fall into a given sub-interval. In such a situation, for example, when two numbers are very close to each other, they can be treated as the same value. In this case, $E = 0$ instead of $E = 1$. In the calculations performed in this study, the variability range of $\alpha_{out}$ was divided into 100 sub-intervals.

Fig. 9 shows a diagram of information entropy as a function of switching time $\tau_r$.

The principle of creating this diagram is the same as creating bifurcation and spectral diagrams. However, neither the time series nor the spectral set are used here, but the entropy value calculated by the formula (13).

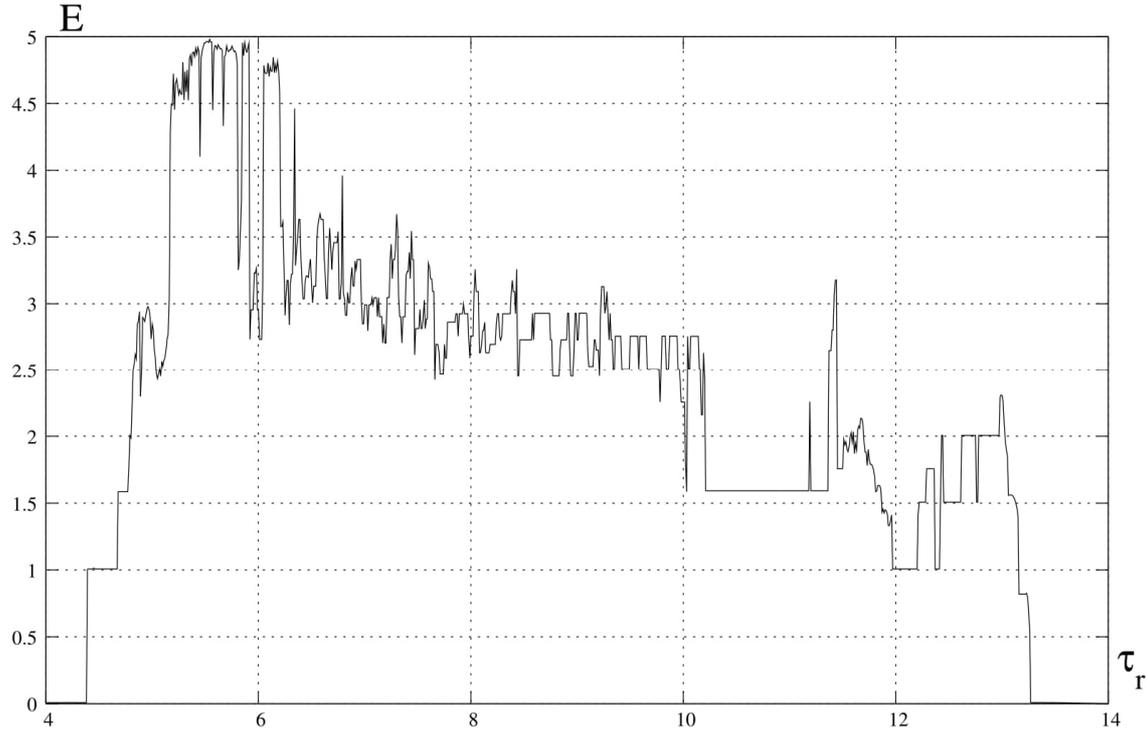

Fig. 9. Diagram of information entropy for $Da = 0.13$

The diagram above shows that the system gives the most information in the chaotic range, i.e. in the range $5.4 < \tau_r < 5.9$.

## 6. Summary

The presented work concerns the analysis of the dynamics of a catalyst (pseudohomogenous) tubular chemical reactor with longitudinal dispersion, in which cyclic transfer of the raw material stream was used. These types of models have been studied in the scientific literature, but neither spectral nor entropy analysis was used there, which was done in this paper. Spectral analysis gives an overview of the participation of individual harmonics in the periodic and chaotic oscillations of the reactor. This information may be useful at the design stage of such a system. In turn, the entropy analysis provides information about the model of the tested system, which may be useful at the stage of modeling this type of apparatus.

The results of the above-mentioned analyzes are presented in phase diagrams, spectral diagrams and in the information entropy diagram. An exemplary Poincare cross-section is also shown, proving the existence of chaotic oscillations in the reactor. Since the reverse system was treated as a discrete system, the above-mentioned cross-section has the character of a Henon attractor. This discretization consists in sampling the values of the reactor variables at the moments of metastasis. In the course of numerical analyzes it turned out that the tested system has the so-called multiple dynamic states. They are characterized by the fact that for the given values of the model parameters, the system can generate various types of oscillations (periodic or non-periodic), depending on the initial conditions of the apparatus. To illustrate this phenomenon, a map of the initial conditions of the reactor was made, on which points determining periodic oscillations were marked. The remaining fields are for chaotic oscillations.

On the basis of the obtained results, three ways of reaching chaos can be noticed. One of them is period doubling, which can be seen on the left side of Fig. 1. Then, around $\tau_r = 5.25$, there is a second way of reaching chaos, namely high amplitude chaos generation, which indicates a crisis. In turn, Fig. 2 shows that for Da = 0.15 there is a third type of reaching chaos, namely through the decay of the torus into the chaotic attractor, i.e. the transition from quasiperiodic solutions to chaotic solutions. This work is purely theoretical. It does not concern the design of a specific apparatus with a specific chemical reaction. It is general in nature, and the adopted parameter values are only examples. As for practical applications with the use of our research, they should take place at the design stage of the reactor system. We have shown: 1. the possibility of certain dynamic phenomena occurring in the presented reactor system, 2: methods of studying dynamic phenomena, including chaotic oscillations.

The use, at the design stage, of the amplitude spectrum method allows, inter alia, to determine the harmonics of the system, i.e. their number, position and amplitude values of these harmonics. This allows, for example, to assess whether the given harmonics pose a process hazard. Such threats may include, for example, harmonics of high frequencies and large amplitudes. The study of the spectrum may also allow, for example, to filter out a given harmonic with a large amplitude, which may be useful especially when one wants to achieve a high average degree of conversion. The study of the entropy of the system allows, in turn, to assess for which values of the system parameters (e.g. for which values of the switching time), this entropy is the highest. Then the system generates the most information about itself, which is included, for example, in the set of measurement results. This maximum information allows for the identification of the model of the tested object on the basis of the above-mentioned results. The analysis of bifurcation diagrams, in turn, allows the evaluation of the values of the system state variables. All this is necessary to design the system in such a way that it works stably and gives the desired process results. From this point of view, this work has unquestionably practical aspects. The "engineering" details were disclosed in Notation. As for the kinetics of the chemical reaction presented in the paper, it should be emphasized that it does not concern a specific process, it is a general exothermic reaction type A→B of any order. If necessary, each time a specific reaction can be substituted for this kinetics and the model can be re-examined.

This work shows that, first of all, complicated dynamic phenomena, including chaotic ones, may occur in the system under study. Secondly, the tendency of these phenomena is clearly shown in the relevant diagrams. For example, the diagram in Fig. 1 clearly shows that chaos occurs within the limits of 4.8 <$\tau_r$= <6.2. For other values of this parameter there are complicated periodic oscillations, and for $\tau_r$= = 13 there are pseudo periodic oscillations. So, the conclusion is that chaos does not exist for large and small switching times! They are somewhere in the middle. This is a very practical point. The diagram in Fig. 2 shows that for 0.14 <$Da$ <0.1675 only periodic or pseudo-periodic oscillations occur. Chaotic oscillations are in the range of 0.075 <*Da* <0.13. For *Da* <0.075 and for *Da*> 0.1675 there is no oscillation at all. So this is another very specific conclusion to be expected when changing the value of *Da*. The same is true of the spectral and entropy diagrams as mentioned above. In turn, Fig. 8 shows for which initial values of state variables the system generates chaotic oscillations, and for which it does not. This means that the tested system has multiple solutions. That is, for the given values of the model parameters, the state of the reactor may be different, depending only on the initial values. Since chaotic oscillations are, as a rule, unfavorable in terms of processes, the research proposed in our article allows to select the system parameters so as to avoid these oscillations in practice. On the contrary, if the designer aims to identify the model, he should select the parameters of the system so that it generates chaotic oscillations with high entropy.

## Notations

| | |
|---|---|
| $c_p$ | heat capacity, $kJ/(kg\ K)$ |
| $C_A$ | concentration of component A, $kmol/m^3$ |
| $Da$ | Damköhler number $\left(=\dfrac{V_R(-r_0)}{\dot{F}C_0}\right)$ |
| $E$ | activation energy, $kJ/kmol$ |
| $\dot{F}$ | volumetric flow rate, $m^3/s$ |
| $(-\Delta H)$ | heat of reaction, $kJ/kmol$ |
| $k$ | reaction rate constant, $1/s(m^3/kmol)^{n-1}$ |
| $L$ | length, $m$ |
| $Le$ | Lewis number, $1+\dfrac{\rho_s c_{ps}}{\rho c_p}$ |
| $m$ | order of reaction |
| $(-r)$ | rate of reaction, $\left(= kC^n\right)$, $kmol/(m^3 s)$ |
| $Pe$ | Peclet number |
| $R$ | gas constant, $kJ/(kmol\ K)$ |
| $t$ | time, $s$ |
| $T$ | temperature, $K$ |
| $V$ | volume, $m^3$ |
| $z$ | position, $m$ |

*Greek letters*

| | |
|---|---|
| $\alpha$ | degree of conversion $\left(= \dfrac{C_{A0} - C_A}{C_{A0}}\right)$ |
| $\beta$ | dimensionless number related to adiabatic temperature increase $\left(= \dfrac{(-\Delta H)C_{A0}}{T_0 \rho c_p}\right)$ |
| $\gamma$ | dimensionless number related to activation energy $\left(= \dfrac{E}{RT_0}\right)$ |
| $\delta$ | dimensionless heat exchange coefficient $\left(= \dfrac{A_q k_q}{\rho c_p \dot{F}}\right)$ |
| $\Theta$ | dimensionless temperature $\left(= \dfrac{T-T_0}{T_0}\right)$ |
| $\xi$ | dimensionless position $\left(= \dfrac{z}{L}\right)$ |

| | | |
|---|---|---|
| $\rho$ | density | $\left(=\dfrac{kg}{m^3}\right)$ |
| $\tau$ | dimensionless time | $\left(=\dfrac{\dot{F}}{V_R}t\right)$ |

*Subscripts*

| | |
|---|---|
| *0* | refers to feed |
| *H* | refers to heat |
| *M* | refers to mass |
| *out* | output of system |
| *r* | refers to reverse flow; switching times |
| *R* | refers to reactor |
| *s* | refers to solid phase |